\documentclass[a4paper,12pt]{article}
\usepackage{fancyhdr,graphicx}
\usepackage{amsfonts}
\usepackage{amssymb}
\usepackage{amsthm}
\usepackage{newlfont}
\usepackage{amsmath}
\usepackage[top=2.0cm,bottom=3.0cm,left=2.5cm,right=2.5cm]{geometry}
\usepackage{multirow}
\usepackage{array}
\usepackage{longtable}
\usepackage{bm}
\usepackage{latexsym}
\usepackage{amsmath}
\usepackage{graphicx}
\usepackage{amssymb}
\usepackage{mathrsfs}
\usepackage{setspace}
\usepackage{amssymb,amsfonts}
\allowdisplaybreaks \numberwithin{equation}{section}

\usepackage{lineno}
\usepackage{algorithm}
\usepackage{algorithmic}
\usepackage{xcolor}
\usepackage{float}
\usepackage{graphics}
\usepackage{graphicx}
\usepackage{epsfig}
\usepackage{booktabs,multirow}
\usepackage{diagbox}

\begin{document}

\title{{
On a conjecture of E\v{g}ecio\v{g}lu and Ir\v{s}i\v{c}
\thanks{This work is supported by Shandong Provincial Natural Science Foundation (ZR2021MA071, ZR2019YQ02)
and National Natural Science Foundation of China (12171414).}}}

\author{
Jianxin Wei $^{a}$\footnote{Corresponding author.},
Yujun Yang $^{b}$\\
\scriptsize{$^{a}$
School of Mathematics and Statistics Science, Ludong University, Yantai, Shandong, 264025, P.R. China}\\
\scriptsize{$^{b}$
School of Mathematics and Information Science, Yantai University, Yantai, Shandong, 264005, P.R. China}\\
\scriptsize{E-mails:
wjx0426@ldu.edu.cn,
yangyj@ytu.edu.cn}}

\date{}

\maketitle

\begin{abstract}
In 2021, {\"O}. E\v{g}ecio\v{g}lu, V. Ir\v{s}i\v{c} introduced the concept of Fibonacci-run graph $\mathcal{R}_{n}$ as an induced subgraph of
Hypercube. They conjectured that the diameter of $\mathcal{R}_{n}$ is given by
$n-\lfloor(1+\frac{n}{2})^{\frac{1}{2}}-\frac{3}{4}\rfloor$.
In this paper, we introduce the novel concept of distance-barriers between vertices in $\mathcal{R}_{n}$ and provide an elegant method to give lower bound for the diameter of $\mathcal{R}_{n}$ via distance-barriers. By constructing different types of distance-barriers, we show that the conjecture does not hold for all $n\geq 230$ and some of $n$ between $91$ and $229$. Furthermore, lower bounds for the diameter of some Fibonacci-run graphs are obtained, which turn out to be better than the result given in the conjecture.

\end{abstract}

\newcommand{\trou}{\vspace{1.5 mm}}
\newcommand{\noi}{\noindent}
\newcommand{\ol}{\overline}
\textbf{Key words:} Hypercube, Fibonacci cube, Fibonacci-run graph, diameter

\section{Introduction}

Let $s$ and $t$ be integers such that $s<t$.
For convenience, we use $s:t$ to denote the set $\{s,s+1,\ldots,t\}$ in this paper.
Let $a$ be any number.
Then we use $\lfloor a\rfloor$ to denote the nearest integer less than or equal to $a$,
and $\lceil a\rceil$ the nearest integer greater than or equal to $a$, respectively.
Let $n\geq1$ and

$\mathcal{B}_{n}=\{b_{1}b_{2}\ldots b_{n} \mid b_{i}\in \{0,1\},i\in 1:n\}$.

\noindent
Then any element of $\mathcal{B}_{n}$ is called
\emph{a binary string} (or simply a string) of length $n$.
The \emph{Hypercube} $Q_{n}$ is the graph defined on
the vertex set $\mathcal{B}_{n}$,
two vertices $b_{1}b_{2}\ldots b_{n}$ and
$b_{1}'b_{2}'\ldots b_{n}'$ being adjacent
if $b_{i}\neq b_{i}'$ holds for exactly one $i$ such that $i\in 1:n$.

A well studied subfamily of hypercubes called Fibonacci cubes were introduced by Hsu \cite{Hsu}.
For $n\geq2$, let

$\mathcal{F}_{n}
=\{b_{1}b_{2}\ldots b_{n}\in\mathcal{B}_{n}\mid b_{i}b_{i+1}=0, 1\leq i\leq n-1\}$,

\noindent
in other words, $\mathcal{F}_{n}$ is the set consisting of all strings of length $n$ that have no two consecutive 1s.
Then an element of $\mathcal{F}_{n}$ is called a \emph{Fibonacci string}.
The \emph{Fibonacci cube} $\Gamma_{n}$ has $\mathcal{F}_{n}$
as vertex set and two vertices are adjacent
if they differ in exactly one coordinate.
For convenience,
we set $\Gamma_{0}\cong K_{1}$ and $\Gamma_{1}\cong K_{2}$.
Fibonacci cubes have been studied from numerous points of view,
as surveyed in \cite{K1}.
Various families of subgraphs of hypercubes have been introduced and studied,
such as in
\cite{EA,ESS,IKR2,IKR,KM,MS, Mu,QW,WYu,Wu,WZ,WY},
and so on.

Fibonacci-run graphs are a new family of graphs introduced by E\v{g}ecio\v{g}lu and Ir\v{s}i\v{c} \cite{EI1,EI2}.
A string $\gamma\in\mathcal{B}_{n}$ is called a \emph{run-constrained string}
if every run of $1$s in $\gamma$ is immediately followed by a strictly longer run of 0s.
Note that a run-constrained string of length $n\geq2$ must be ended with $00$.
For $n\geq0$,
the \emph{Fibonacci-run graph} $\mathcal{R}_{n}$
has the vertex set

$V(\mathcal{R}_{n})=\{r\mid r00$
is a run-constrained binary string of length $n+2\}$,

\noindent
and two vertices are adjacent
if they differ in exactly one coordinate.
Note that if $n=0$,
then $V(\mathcal{R}_{n})= \{\lambda\}$ and so $\mathcal{R}_{n}\cong K_{1}$,
where $\lambda$ is the null string.
If $n=1$,
then it is easily seen that $\mathcal{R}_{1}\cong K_{2}$.

Although the Fibonacci-run graph $\mathcal{R}_{n}$ has the
same number of vertices as the Fibonacci cube $\Gamma_{n}$,
it has fewer edges than $\Gamma_{n}$, as well as some other different graph
theoretical properties from $\Gamma_{n}$, as shown in papers \cite{EI1,EI2}.
As the study on $\mathcal{R}_{n}$ is in the initial stage,
various crucial graph properties for $\mathcal{R}_{n}$ need to be determined or characterized.
In \cite{EI1,EI2},
some interesting questions and conjectures on Fibonacci-run graphs have been proposed.
For example,
what is the number of vertices of degree $k$ in $\mathcal{R}_{n}$ and
what is the irregularity of $\mathcal{R}_{n}$?
As for the hamiltonicity,
it was conjectured that $\mathcal{R}_{n}$ is Hamiltonian if and only if $n\equiv1$ (mod 3).
As for the distance parameters of graph,
they asked what is the radius Fibonacci-run graph $\mathcal{R}_{n}$ and conjectured
and the diameter of $\mathcal{R}_{n}$ is given by

\trou \noi {\bf Conjecture 1.1\cite{EI1}.}
\emph{Let $n\geq0$. Then diam$(\mathcal{R}_{n})=n-\lfloor\sqrt{1+\frac{n}{2}}-\frac{3}{4}\rfloor$.}
\\The problem on the radius of $\mathcal{R}_{n}$ has been completely solved in \cite{W0}. In the present paper, we focus on conjecture 1.1.
It is shown that the above conjecture is not true for all the dimensions $n\in \{91,94,95,119:124,131,136,152:170,172,178,181,184,185,189:223\}$ and $n\geq 230$ in section 4.
The main tool used is the so called distance-barriers between vertices of $\mathcal{R}_{n}$,
as introduced in Section 3.

In the rest of this section,
we give some necessary notations and concepts.
Let $\alpha\in\mathcal{B}_{n}$ be any string.
Then we use $|\alpha|$ to denote the length of $\alpha$.
For convenience, we use the product to denote string concatenation,
for example, $1^{s}$ is the string $11\ldots1$ of length $s$.
More generally, for a string $\alpha$,
$\alpha^{s}$ is the concatenation of $s$ copies of $\alpha$.
We also take the null string $\lambda$ into consideration.
For examples,
$1^{0}=\lambda$ and in generally,
$\alpha^{0}$ also is null string.
A non-extendable sequence of contiguous equal bits in a string is called a \emph{run} (or \emph{block}) of the string.
For example,
$1^{3}0^{2}1^{2}0$ has 4 runs.
A string $f$ is called a \emph{factor} of a string $\alpha$
if $f$ appears as a sequence of $|f|$ consecutive bits of $\alpha$.

\section{Preliminaries}

In this section, we introduce some distance-based concepts.
Note that the \emph{distance} $d_{G}(u,v)$ (or simply $d(u,v)$)
between vertices $u$ and $v$ of a graph $G$ is the length of
a shortest path between them in $G$.

The \emph{eccentricity} $e(v)$ of a vertex $v$ in a connected graph $G$
is the maximum distance from $v$ and any other vertex, i.e.

$e(v)=\max\limits_{u\in V(G)} d(v,u)$.

The \emph{radius} $rad(G)$ of $G$ is
the minimum eccentricity of the vertices of $G$,
i.e.

rad$(G)=\min\limits_{v\in V(G)} e(v)$.

\noindent
A vertex is called \emph{central} if $e(v)=$ rad$(G)$.
The \emph{center} $Z(G)$ of $G$ is the set of all central vertices.

The radius and center of Fibonacci cube $\Gamma_{n}$ and Fibonacci-run graph $\mathcal{R}_{n}$ are determined in \cite{MS} and \cite{W0},
respectively.
It is interesting to note that for a given $n\geq0$,
$\mathcal{R}_{n}$ and $\Gamma_{n}$ have the same radius $\lceil\frac{n}{2}\rceil$.

The \emph{diameter} $diam(G)$ of $G$ is the maximal eccentricity $e(v)$ when $v$ runs in $G$, i.e.

$diam(G)=\max\limits_{v\in V(G)} e(v)=\max\limits_{u,v\in V(G)} d(u,v)$.


Hsu \cite{Hsu} proved that the diameter of $\Gamma_{n}$ is $n$.
E\v{g}ecio\v{g}lu and Ir\v{s}i\v{c} \cite{EI1} conjectured the diameter of $\mathcal{R}_{n}$ is given by $n-\lfloor\sqrt{1+\frac{n}{2}}-\frac{3}{4}\rfloor$ as shown in Conjecture 1.1.

The \emph{Hamming distance} $H(\alpha,\beta)$ between two strings
$\alpha,\beta\in\mathcal{B}_{n}$ is the number of coordinates
in which $\alpha$ and $\beta$ differ.
It is well known that $d_{Q_{n}}(\alpha,\beta)=H(\alpha,\beta)$
for any vertices $\alpha,\beta \in V(Q_{n})$ \cite{HIK}.

Let $H$ and $G$ be connected graphs.
Then a mapping $f:V(H)\rightarrow V(G)$ is an \emph{isometric embedding} if $d_{H}(u,v)=d_{G}(f(u),f(v))$ for any $u,v\in V(H)$.
A \emph{partial cube} is a connected graph which admits an isometric embedding into a hypercube \cite{HIK}.
It is clear that the distance between any two vertices of a partial cube is exactly the Hamming distance between them.

It is shown that $\Gamma_{n}$ is an isometric subgraph of $Q_{n}$ \cite{Hsu},
that is,
$\Gamma_{n}$ is a partial cube.
However,
Fibonacci-run graphs are in general not partial cubes \cite{EI1},
although there exist some particular vertices $\alpha,\beta$ of $\mathcal{R}_{n}$ with $d_{\mathcal{R}_{n}}(\alpha,\beta)=H(\alpha,\beta)$.
In the following, we use an example to illustrate that there exists $\alpha,\beta$ of $\mathcal{R}_{n}$ with $d_{\mathcal{R}_{n}}(\alpha,\beta)=H(\alpha,\beta)$.

\trou \noi {\bf Example 2.1.}
Let $\alpha$ and $\beta$ be two vertices of $\mathcal{R}_{28}$ such that
$$\alpha=1001111111000000001100011100 \mbox{~~and~~} \beta=1110000100111000011111100000.$$
Then $H(\alpha,\beta)=18$ and we can find a $\alpha,\beta$-path of length 18 in $\mathcal{R}_{n}$ as follows:

$1001111111000000001100011100\rightarrow$
$1000111111000000001100011100\rightarrow$

$1000011111000000001100011100\rightarrow$
$1000001111000000001100011100\rightarrow$

$1000000111000000001100011100\rightarrow$
$1000000110000000001100011100\rightarrow$

$1000000100000000001100011100\rightarrow$
$1000000100000000001100001100\rightarrow$

$1000000100000000001100000100\rightarrow$
$1000000100000000001100000000\rightarrow$

$1100000100000000001100000000\rightarrow$
$1110000100000000001100000000\rightarrow$

$1110000100100000001100000000\rightarrow$
$1110000100110000001100000000\rightarrow$

$1110000100111000001100000000\rightarrow$
$1110000100111000011100000000\rightarrow$

$1110000100111000011110000000\rightarrow$
$1110000100111000011111000000\rightarrow$

$1110000100111000011111100000$.

From the above example we can see that such a shortest path is not a common \emph{canonical} $\alpha,\beta$-path
(that is, a path where we first change each bit of $\alpha$ from $1$ to $0$ for which $a_{i}=1$ and $b_{i}=0$,
and then change from $0$ to $1$ the other bits in which $\alpha$ and $\beta$ differ)
but a path obtained by using a particular order to change the bits in which they differ.
However,
a shortest path (lying completely in $\mathcal{R}_{n}$) with Hamming distance does not always exist between some vertices of $\mathcal{R}_{n}$ (as shown in Example 3.1).
This will be discussed in the following section.

\section{Distance-barrier of $\mathcal{R}_{n}$}

This section contains two subsections. In the first subsection,
we introduce the concept of distance-barrier between two vertices of $\mathcal{R}_{n}$,
which will be used to give sufficient and necessary conditions for the existence of the shortest path with Hamming distance between vertices of $\mathcal{R}_{n}$.
In the second subsection,
we apply those results to study $diam(\mathcal{R}_{n})$.

\subsection{Basic properties of distance-barriers}

First we give an example to show that there is not always a shortest path with Hamming distance between two vertices in $\mathcal{R}_{n}$.

\trou \noi {\bf Example 3.1.}
We choose the following two vertices $\beta$ and $\gamma$ of $\mathcal{R}_{21}$:

$\beta=b_{1}b_{2}\ldots b_{20}b_{21}=100110001110000011100$ and

$\gamma=c_{1}c_{2}\ldots c_{20}c_{21}=111111111110000000000$.

\noindent
It is clear that $H(\beta,\gamma)=8$.
Let $A=b_{1}=1$,
$B=b_{4}b_{5}=11$ and $C=b_{9}b_{10}b_{11}=111$ be the first three factors consisting of $1$s of $\beta$.
Then it can be seen that if we want to get a path between $\beta$ and $\gamma$ lying completely in $\mathcal{R}_{21}$,
we must change all the bits contained in at least two of $A,B$ and $C$ from $1$ to $0$ and then from $0$ to $1$.
Obviously,
by the fact that $|A|<|B|<|C|$,
a shortest path between $\beta$ and $\gamma$ (lying in $\mathcal{R}_{21}$) could be obtained by changing twice times the bits in $A$ and $B$ as follows:

$100110001110000011100\rightarrow$
$000110001110000011100\rightarrow$
$000010001110000011100\rightarrow$

$000000001110000011100\rightarrow$
$000000011110000011100\rightarrow$
$000000111110000011100\rightarrow$

$000001111110000011100\rightarrow$
$000011111110000011100\rightarrow$
$000111111110000011100\rightarrow$

$001111111110000011100\rightarrow$
$011111111110000011100\rightarrow$
$111111111110000011100\rightarrow$

$111111111110000001100\rightarrow$
$111111111110000000100\rightarrow$
$111111111110000000000$.

\noindent
It is easily seen that the distance between $\beta$ and $\gamma$ is

$d_{\mathcal{R}_{21}}(\beta,\gamma)$

$=H(\beta,\gamma)+2(|A|+|B|)$

$=8+2(1+2)$

$=14.$

\noindent
From the above equations we know that $d_{\mathcal{R}_{21}}(\beta,\gamma)=14$, which is not equal to $H(\beta,\gamma)=8$.

According to the characters of the factors
$\beta'=b_{1}b_{2}\ldots b_{11}=10011000111$ of $\beta$ and
$\gamma'=c_{1}c_{2}\ldots c_{11}=11111111111$ of $\gamma$ in Example 3.1,
we are motivated to define the following useful concept.

\trou \noi {\bf Definition 3.2.}
\emph{Let $\alpha=a_{1}a_{2}\ldots a_{n}$ and $\beta=b_{1}b_{2}\ldots b_{n}$ be two vertices of $\mathcal{R}_{n}$.
Suppose that for some $t\geq1$ there exist factors
$$\alpha'=1^{r_{1}}0^{s_{1}}1^{r_{2}}\ldots0^{s_{t}}1^{r_{t+1}} \mbox{~~and~~} \beta'=1^{r_{1}}1^{s_{1}}1^{r_{2}}\ldots1^{s_{t}}1^{r_{t+1}}$$
of $\alpha$ and $\beta$ respectively such that the $k$-th coordinate of $\alpha$ (resp. $\beta$) is the first coordinate from which the factor $\alpha'$ (resp. $\beta'$) begins.
If $a_{k-1}+b_{k-1}=1$ (or $k-1=0$) and $a_{k+\mid\alpha'\mid}+b_{k+\mid\alpha'\mid}=1$ (or $k+\mid\alpha'\mid>n$),
then we call
$\begin{bmatrix}
\alpha'\\
\beta'
\end{bmatrix}$
a \emph{distance-barrier of length $|\alpha'|$} between $\alpha$ and $\beta$ with $t+1$ barriers. In addition, we call
$1^{r_{i}}$ the $i$-th \emph{barrier} with
$r_{i}$ being the \emph{thickness} of barrier $1^{r_{i}}$,
$i\in1:(t+1)$.}

For example, by Definition 3.2,
we know that
$\begin{bmatrix}
\beta'\\
\gamma'
\end{bmatrix}=$
$\begin{bmatrix}
10011000111\\
11111111111
\end{bmatrix}$
is a distance-barrier of length $11$ between vertices $\beta$ and $\gamma$ (in Example 3.1) of $\mathcal{R}_{11}$ with three barriers of thickness $1$, $2$ and $3$, respectively.

Now we consider the effect of a distance-barrier between two vertices on the distance between them in the Fibonacci-run graph.

Suppose that
$\begin{bmatrix}
\alpha'\\
\beta'
\end{bmatrix}$
$=\begin{bmatrix}
1^{r_{1}}0^{s_{1}}1^{r_{2}}\ldots0^{s_{t}}1^{r_{t+1}}\\
1^{r_{1}}1^{s_{1}}1^{r_{2}}\ldots1^{s_{t}}1^{r_{t+1}}
\end{bmatrix}$
is a distance-barrier between vertices $\alpha$ and $\beta$ of $\mathcal{R}_{n}$,
and $M_{(\alpha',\beta')}=\max\{r_{1},\ldots,r_{t},r_{t+1}\}$.
As $t$ barriers must be changed at least twice to get a path lying $\mathcal{R}_{n}$ between $\alpha$ and $\beta$,
we should change the $t$ barriers with the smallest thickness twice to get a shortest path between $\alpha$ and $\beta$.
Let $C_{(\alpha',\beta')}
=\mathop{\sum}\limits_{i\in 1:t+1\setminus M_{(\alpha',\beta')}}r_{i}$
and $C'_{(\alpha',\beta')}
=\mathop{\sum}\limits_{i\in 1:t+1}r_{i}$.
Obviously,
$C'_{(\alpha',\beta')}-C_{(\alpha',\beta')}=M_{(\alpha',\beta')}$.
Then we have the following result which are used to determine the length of a shortest path between two vertices if there are distance-barriers between them.

\trou \noi {\bf Lemma 3.3.}
\emph{Let $\begin{bmatrix}
\alpha'_{k}\\
\beta'_{k}
\end{bmatrix}$ be a distance-barrier between vertices $\alpha$ and $\beta$ of $\mathcal{R}_{n}$,
$k\in1:s$.
Then $d_{\mathcal{R}_{n}}(\alpha,\beta)
=H(\alpha,\beta)+2\mathop{\sum}\limits^{s}\limits_{k=1}C_{(\alpha'_{k},\beta'_{k})}
=H(\alpha,\beta)+2\mathop{\sum}\limits^{s}\limits_{k=1}(C'_{(\alpha'_{k},\beta'_{k})}
-M_{(\alpha'_{k},\beta'_{k})})$.}

The following result is the main tool used in the present paper.

\trou \noi {\bf Lemma 3.4.}
\emph{Let $\alpha$ and $\beta$ be two vertices of $\mathcal{R}_{n}$.
Then $d_{\mathcal{R}_{n}}(\alpha,\beta)=H(\alpha,\beta)$ if and only if there is no distance-barrier between $\alpha$ and $\beta$.}

\trou \noi {\bf Proof.}
If there is a distance-barrier between $\alpha$ and $\beta$,
then obviously $d_{\mathcal{R}_{n}}(\alpha,\beta)>H(\alpha,\beta)$ by Lemma 3.3.

Now we suppose that there is no distance-barrier between $\alpha=a_{1}a_{2}\ldots a_{n}$ and $\beta=b_{1}b_{2}\ldots b_{n}\in V(\mathcal{R}_{n})$.
Suppose that $H(\alpha,\beta)=h$ and $a_{i_{j}}+b_{i_{j}}=1$ for $j\in1:h$.
Now we show that there is a $\alpha,\beta$-path with Hamming distance $H(\alpha,\beta)$ lying completely in $\mathcal{R}_{n}$ by induction on $H(\alpha,\beta)$.

Obviously,
the assertion holds for $H(\alpha,\beta)=1$.
Now we suppose that it holds for vertices with Hamming distance $h-1$ ($h\geq2$),
that is, if there is no distance-barrier between vertices with Hamming distance $h-1$,
then there is a path with distance $h-1$ lying in $\mathcal{R}_{n}$.
Without loss of generality,
suppose that $a_{i_{1}}=1$ and so $b_{i_{1}}=0$.
For convenience, we distinguish the following three cases:

$(1)$ $i_{1}=1$ or $i_{1}=n$.

$(2)$ $a_{i_{1}-1}=0$ or $a_{i_{1}+1}=0$.

$(3)$ $a_{i_{1}-1}=1$ and $a_{i_{1}+1}=1$

For Cases $(1)$ and $(2)$,
let $\alpha'$ be the string obtained by changing the $i_{1}$-th bit of $\alpha$ from $0$ to $1$.
For Case $(3)$, as $a_{i_{1}}=1$,
it is known that $b_{i_{1}-1}=1$.
Since there is no distance-barrier between vertices $\alpha$ and $\beta$,
there must exist $s\geq1$ such that $a_{i_{1}+1}=\ldots=a_{i_{1}+s}=1$ and $b_{i_{1}+1}=\ldots=b_{i_{1}+s}=0$,
but $a_{i_{1}+s+1}=0$.
For this case,
let $\alpha'$ be the string obtained by changing the $(i_{1}+s)$-th bit of $\alpha$ from $0$ to $1$.

It is clear that the string $\alpha'$ obtained in the above three cases is a vertex of $\mathcal{R}_{n}$ with
$H(\alpha',\beta)=h-1$. Thus there is no distance-barrier between $\alpha'$ and $\beta$.
Therefore,
there is a shortest $\alpha',\beta$-path with distance $h-1$ lying in $\mathcal{R}_{n}$.
Hence we know that if there is no distance barrier between $\alpha$ and $\beta$,
then $d_{\mathcal{R}_{n}}(\alpha,\beta)
=1+d_{\mathcal{R}_{n}}(\alpha',\beta)
=1+H(\alpha',\beta)
=H(\alpha,\beta)$.
This completes the proof.
$\Box$

\trou \noi {\bf Corollary: 3.5\cite{EI1}.}
\emph{$\mathcal{R}_{n}$ is a partial cube if and only if $n\leq 6$.}

\trou \noi {\bf Proof.}
It is obvious that the distance-barrier with the smallest length is
$\begin{bmatrix}
\alpha'\\
\beta'
\end{bmatrix}$
$=\begin{bmatrix}
1001\\
1111
\end{bmatrix}$,
so the length of the run-constrained string contain $1111$ as a factor must be at least $7$.
If $n\geq7$,
then $\alpha=10010000^{n-7}$ and $\beta=11110000^{n-7}$ are vertices of $\mathcal{R}_{n}$,
and so
$\begin{bmatrix}
\alpha'\\
\beta'
\end{bmatrix}$
is a distance-barrier between $\alpha$ and $\beta$.
The desired result follows by Lemma 3.4.
$\Box$

\subsection{Applications of distance-barrier}

A lower bound for the diameter of $\mathcal{R}_{n}$ for some $n$ is given by E\v{g}ecio\v{g}lu and Ir\v{s}i\v{c} (\cite{EI1}, Lemma 5.1 and Theorem 5.3).
On the basis of this result and some calculation examples,
they proposed Conjecture 1.1.
In this section,
firstly, a general lower bound for the diameter of $\mathcal{R}_{n}$ for all $n\in\mathbb{N}$ is obtained (Lemma 3.6).
Then by constructing vertices with distance-barrier between them in various ways, we obtain some improved lower bounds for the diameter of $\mathcal{R}_{n}$ (see Lemmas 3.8-3.10).

\trou \noi {\bf Lemma 3.6.}
\emph{Let $n$ be any non-negative integer.
If $n\in (2p^{2}+3p):(2p^{2}+7p+4)$ for some $p\geq0$,
then $diam(\mathcal{R}_{n})\geq n-p$.}

\trou \noi {\bf Proof.}
Let $p\geq0$ be an integer,
$n=2p^{2}+7p+4$,
and $\mu$ and $\nu$ be the following two vertices of $\mathcal{R}_{n}$:

$\mu=10^{2}1^{3}\ldots1^{2p+1}0^{2p+2}1^{p+1}0^{p}$,

$\nu=01^{2}0^{3}\ldots0^{2p+1}1^{2p+2}0^{p+1}0^{p}$.
~~~~~~~~~~~~~~~~~~~~~~~~~~~~~~~~~~~~~~~~~~~~~~~~~~~~~~~~~~~~~~~~~~(3.1)

\noindent
Note that if $p=0$,
then

$\mu=1001$,

$\nu=0110$;

\noindent
and if $p=1$,
then

$\mu=1001110000110$,

$\nu=0110001111000$.

As there is no distance-barrier between $\mu$ and $\nu$,
$diam(\mathcal{R}_{n})\geq d_{\mathcal{R}_{n}}(\mu,\nu)=H(\mu,\nu)=n-p$ by Lemma 3.4.
Since any suffix of $\mu$ or $\nu$ is also a run-constrained string,
we know that for any $n\leq2p^{2}+7p+4$, $diam(\mathcal{R}_{n})\geq n-p$.
By the fact that $2(p-1)^{2}+7(p-1)+4=2p^{2}+3p-1$,
we know that the set of non-negative integers
$\mathbb{N}=\bigcup\limits^{+\infty}\limits_{p=0}(2p^{2}+3p):(2p^{2}+7p+4)$.
The result holds.
$\Box$

For convenience,
the pair of vertices $\mu$ and $\nu$ (in Eq. (3.1)) or their suffixes of the same length $n(\leq2p^{2}+7p+4)$ in the proof of Lemma 3.6 are called \emph{a pair of $H$-type strings ending with} $p$ \emph{of length} $n$.
The vertices $\mu$ and $\nu$ (of length $2p^{2}+7p+4$) are called a pair of \emph{complete} $H$-type strings ending with $p$.
Otherwise,
the $H$-type strings of length less than $2p^{2}+7p+4$ is called \emph{incomplete}.
It is easily seen that the length of the longest run of the $H$-type string ending with $p$ is $2(p+1)$.
As there is no distance-barrier between a pair of $H$-type strings,
the distance beween a pair of $H$-type strings of length $n$ in $\mathcal{R}_{n}$ is $n-p$ by Lemma 3.4.

The lower bound given in Lemma 3.6 is obtained by the distance between two vertices of $\mathcal{R}_{n}$ with no distance-barrier between them.
By Lemma 3.3,
we have every reason to believe that the bound could be improved.
In the following, we first give an example to indicate that distance-barrier can be used as a tool in studying the diameter of $\mathcal{R}_{n}$.
Then we give a specific type of distance-barriers which will be used in the next section.

\trou \noi {\bf Example 3.7.}
Let $n=91=2p^{2}+3p+1$,
where $p=6$.
Then by Lemma 3.6 we know for the following incomplete $H$-type strings ending with $p$

$\mu=10^{8}1^{9}0^{10}1^{11}0^{12}1^{13}0^{14}1^{7}0^{6}$ and

$\nu=01^{8}0^{9}1^{10}0^{11}1^{12}0^{13}0^{14}0^{7}0^{6}$,

\noindent
$diam(\mathcal{R}_{n})\geq d_{\mathcal{R}_{n}}(\mu,\nu)=n-6=85.$
If we chose vertices

$\alpha=(111)^{5}10^{17}1^{19}0^{20}1^{10}0^{9}$ and

$\beta=(100)^{5}11^{17}0^{19}1^{20}0^{10}0^{9}$,

\noindent
then there is a distance-barrier
$\begin{bmatrix}
\alpha'\\
\beta'
\end{bmatrix}$
=$\begin{bmatrix}
(111)^{5}1\\
(100)^{5}1
\end{bmatrix}$
between $\alpha$ and $\beta$,
and so $diam(\mathcal{R}_{n})\geq d_{\mathcal{R}_{n}}(\mu,\nu)
=H(\alpha,\beta)+C_{(\alpha',\beta')}
=86$ by Lemma 3.3.
The example ends here.

It can be seen that vertices $\alpha$ and $\beta$ in Example 3.7 could be obtained from the following incomplete $H$-type strings,
respectively:

$1^{15}0^{18}1^{19}0^{20}1^{10}0^{9}$ and

$0^{15}1^{18}0^{19}1^{20}0^{10}0^{9}$.

Inspired by this example,
we give a type of distance-barriers between vertices in $\mathcal{R}_{n}$
which can be used to improve the lower bound for the diameter of $\mathcal{R}_{n}$,
as given in the following result.

\trou \noi {\bf Lemma 3.8.}
\emph{Let $t\geq2$, $s\in0:2$, $k\leq3t+s$, $2\leq b\leq t$ and}

$\mu=0^{k}1^{3t+s+1}0^{3t+s+2}\cdots1^{3t+s+i-1}0^{3t+s+i}\cdots$,

$\nu=1^{k}0^{3t+s+1}1^{3t+s+2}\cdots0^{3t+s+i-1}1^{3t+s+i}\cdots$

\noindent
\emph{be a pair of incomplete H-type strings,
where $i\geq2$ and the length of the longest run contained in $\mu$ (or $\nu$) is at least $3t+s+\lfloor\frac{3b-k-4}{3(t-b)+s+1}\rfloor$.
Then we can construct a pair of vertices such that there is a distance-barrier
$Bar_{1}=$
$\begin{bmatrix}
(111)^{b-1}1\\
(100)^{b-1}1
\end{bmatrix}$
between them.}

\trou \noi {\bf Proof.}
It can be seen that if $3b\leq k\leq3t+s$,
then for the following strings

$\alpha=0^{k-3b}(100)^{b}1^{3t+s+1}0^{3t+s+2}\cdots$ and

$\beta=1^{k-3b}(111)^{b}0^{3t+s+1}1^{3t+s+2}\cdots$,

\noindent
$Bar_{1}$ is a distance-barrier between $\alpha$ and $\beta$.

If $k=3b-1$,
then $Bar_{1}$
is a distance-barrier between strings

$\alpha=11000(100)^{b-2}11^{3t+s-1}0^{3t+s+2}\cdots$ and

$\beta=01111(111)^{b-2}10^{3t+s-1}1^{3t+s+2}\cdots$.

If $3b-4< k<3b-1$,
then
$Bar_{1}$
is a distance-barrier between strings

$\alpha=0^{k-3b+3}(100)^{b-1}11^{3t+s}0^{3t+s+2}\cdots$ and

$\beta=1^{k-3b+3}(111)^{b-1}10^{3t+s}1^{3t+s+2}\cdots$.

For the case that $1\leq k\leq 3b-4$, we first give a disjoint partition of $1:3b-4$.
Let $m=\lfloor\frac{3b-4}{3(t-b)+s+1}\rfloor$
and $I_{i}=(3b-4-i-(i+1)(3(t-b)+s))):(3b-4-i-i(3(t-b)+s)))$,
where $i=0,1,\ldots,m$.
In other words,
starting from $3b-4$, we divide $1:3b-4$ into $m+1$ subsets such that every subset consists of $3(t-b)+s$ consecutive numbers.
Note that when $i=m$,
$I_{i}=1:(3b-4-m-m(3(t-b)+s)))$.
Then
$1:3b-4=\bigcup\limits^{m}\limits_{i=0}I_{i}$.

For $k\in I_{0}$,
that is, $3b-4-[3(t-b)+s]\leq k\leq 3b-4$, let

$\alpha=(100)^{b-1}11^{3(t-b)+k+s+3}0^{3t+s+2}\cdots$,

$\beta=(111)^{b-1}10^{3(t-b)+k+s+3}1^{3t+s+2}\cdots$

\noindent
Then $Bar_{1}$ is a distance-barrier between $\alpha$ and $\beta$.

If $k\in I_{i}$,
$i\in1:m$,
without loss of generality, suppose that $i$ is odd.
Let

$\alpha=(100)^{b-1}11^{3b-1}0^{3b+1}\cdots1^{3b+i-1}0^{(i+1)(3(t-b)+s)+2i+k+4}\cdots$,

$\beta=(111)^{b-1}10^{3b-1}1^{3b+1}\cdots0^{3b+i-1}1^{(i+1)(3(t-b)+s)+2i+k+4}\cdots$.

\noindent
Then $Bar_{1}$ is a distance-barrier between $\alpha$ and $\beta$.
This completes the proof.
$\Box$

Lemma 3.8 implies that based on some pair of $H$-type strings of length long enough
(note that it need $\lfloor\frac{3b-k-4}{3(t-b)+s+1}\rfloor$ runs from the second run to the longest run),
the strings having distance-barrier with $b$ barriers can be constructed.
For the case $s=1,k=3t$ and $b=t$,
we can construct strings having distance-barrier with $b+1$ barriers as shown in the following lemma.

\trou \noi {\bf Lemma 3.9.}
\emph{Let $t\geq2$ and}

$\mu=0^{3t}1^{3t+3}0^{3t+4}\cdots$, and

$\nu=1^{3t}0^{3t+3}1^{3t+4}\cdots$

\noindent
\emph{be a pair of incomplete H-type strings.
Then a pair of vertices with distance-barrier
$Bar_{2}=$
$\begin{bmatrix}
(111)^{t}1\\
(100)^{t}1
\end{bmatrix}$
between them can be constructed.}

\trou \noi {\bf Proof.}
Let

$\alpha=(100)^{t}11^{3t+2}0^{3t+4}\cdots$ and

$\beta=(111)^{t}10^{3t+2}1^{3t+4}\cdots$.

\noindent
Then $Bar_{2}$ is a distance-barrier between $\alpha$ and $\beta$.
$\Box$

Suppose that for some $p\geq0$,
the length of the pair of incomplete $H$-type strings $\mu,\nu$ in Lemmas 3.8 (or Lemma 3.9) is $n$,
where $2p^2+3p\leq n\leq2p^2+7p+4$.
Then we can assume that $\mu$ and $\nu$ end with $p+m$ for some $m\geq0$.
So $d_{\mathcal{R}_{n}}(\mu,\nu)=H(\mu,\nu)=n-(p+m)$ by Lemma 3.4.
For convenience,
vertices $\alpha$ and $\beta$ constructed on the basis of $\mu$ and $\nu$ in the proof of Lemma 3.8 (or Lemma 3.9)
are called \emph{NH-type strings ending with} $p+m$.
Let $\begin{bmatrix}
\alpha'\\
\beta'
\end{bmatrix}$ be the distance-barrier
$\begin{bmatrix}
(111)^{r-1}1\\
(100)^{r-1}1
\end{bmatrix}$
between $\alpha$ and $\beta$ ($r=b$ in Lemma 3.8 and $r=b+1$ in Lemma 3.9).
Then there are $r$ barriers in
$\begin{bmatrix}
\alpha'\\
\beta'
\end{bmatrix}$
and the thickness of every barrier is $1$.
It is easily seen that $H(\mu,\nu)=H(\alpha,\beta)+r$,
and $C_{(\alpha',\beta')}=r-1$.
Therefore,
we could establish the relation  between $d_{\mathcal{R}_{n}}(\alpha,\beta)$ and $d_{\mathcal{R}_{n}}(\mu,\nu)$  by Lemma 3.3:

$d_{\mathcal{R}_{n}}(\alpha,\beta)$

$=H(\alpha,\beta)+2C_{(\alpha',\beta')}$

$=H(\mu,\nu)-r+2(r-1)$

$=d_{\mathcal{R}_{n}}(\mu,\nu)+r-2$

$=n-(p+m)+r-2$

$=(n-p)+r-(m+2)$.~~~~~~~~~~~~~~~~~~~~~~~~~~~~~~~~~~~~~~~~~~~~~~~~~~~~~~~~~~~~~~~~~~~~~~~~~~~~~~~~~~~~~~~~~(3.2)

More generally,
we get the following result.

\trou \noi {\bf Lemma 3.10.}
\emph{Let $\mu$ and $\nu$ be a pair of H-type strings,
and $\alpha$ and $\beta$ be a pair of NH-type strings obtained from $\mu$ and $\nu$.
Suppose $\begin{bmatrix}
\alpha'_{k}\\
\beta'_{k}
\end{bmatrix}$ is the distance-barriers between vertices $\alpha$ and $\beta$,
$k=1,\ldots,s$.
Then}

$d_{\mathcal{R}_{n}}(\alpha,\beta)$

$=d_{\mathcal{R}_{n}}(\mu,\nu)+\mathop{\sum}\limits^{s}\limits_{k=1}(C_{(\alpha'_{k},\beta'_{k})}-M_{(\alpha'_{k},\beta'_{k})})$

$=d_{\mathcal{R}_{n}}(\mu,\nu)+\mathop{\sum}\limits^{s}\limits_{k=1}(C'_{(\alpha'_{k},\beta'_{k})}
-2M_{(\alpha'_{k},\beta'_{k})})$.

\trou \noi {\bf Proof.}
As the barrier between of $\alpha$ and $\beta$ are obtained from $\mu$ and $\nu$ by changing some bits of $\mu$ or $\nu$ from $0$ to $1$,
$H(\mu,\nu)=H(\alpha,\beta)+\mathop{\sum}\limits^{s}\limits_{k=1}C'_{(\alpha'_{k},\beta'_{k})}.$
By Lemma 3.3,
$d_{\mathcal{R}_{n}}(\mu,\nu)=H(\mu,\nu)$ and $C'_{(\alpha'_{k},\beta'_{k})}-C_{(\alpha'_{k},\beta'_{k})}=M_{(\alpha'_{k},\beta'_{k})}$.
So

$d_{\mathcal{R}_{n}}(\alpha,\beta)$

$=H(\alpha,\beta)+2\mathop{\sum}\limits^{s}\limits_{k=1}C_{(\alpha'_{k},\beta'_{k})}$

$=d_{\mathcal{R}_{n}}(\mu,\nu)-\mathop{\sum}\limits^{s}\limits_{k=1}C'_{(\alpha'_{k},\beta'_{k})}
+2\mathop{\sum}\limits^{s}\limits_{k=1}C_{(\alpha'_{k},\beta'_{k})}$,

\noindent
and so the result holds.
$\Box$

Note that Eq. (3.2) is the case of Lemma 3.10 with only one distance-barrier ($s=1$) having $b$ barriers of thickness of 1 between $\alpha$ and $\beta$.

\section{On the diameter of $\mathcal{R}_{n}$}

In this section,
first we show Conjecture 1.1 does not hold for almost all $n\in\mathbb{N}$.
Then some exploring work on the diameter of $\mathcal{R}_{n}$ is carried out.

\subsection{Conjecture 1.1 is in general not true}

In this subsection,
it is shown that when dimension $n\geq91$, the result in Conjecture 1.1 comes to fail for some $n$,
and for all $n\geq230$,
Conjecture 1.1 does not hold.

First,
we show that diameter of $\mathcal{R}_{n}$ given in Conjecture 1.1 is nothing but the lower bound in Lemma 3.6,
that is,
the following result holds.

\trou \noi {\bf Lemma 4.1.}
\emph{Let $n\in (2p^{2}+3p):(2p^{2}+7p+4)$ for some $p\geq0$.
Then $\lfloor\sqrt{1+\frac{n}{2}}-\frac{3}{4}\rfloor$=$p$.}

\trou \noi {\bf Proof.}
Let $a>0$ and $b\geq0$.
Then $\sqrt{a^{2}+b}-a\geq0$ and
it is easily seen that the following result holds:

$\sqrt{a^{2}+b}-a
\begin{cases}
>1,b>2a+1,\\
=1,b=2a+1,\\
<1,b<2a+1.\\
\end{cases}$

Let $n\in (2p^{2}+3p):(2p^{2}+7p+4)$.
Then $n=2p^{2}+3p+k$ for some $k$ such that $0\leq k\leq 4p+4$.
Let $a=p+\frac{3}{4}$ and $b=\frac{k}{2}+\frac{7}{16}$.
Then $b<p+1+\frac{7}{16}<2p+\frac{3}{2}+1=2a+1$.
As

$\sqrt{1+\frac{n}{2}}-\frac{3}{4}-p$

$=\sqrt{1+\frac{2p^{2}+3p+k}{2}}-(\frac{3}{4}+p)$

$=\sqrt{(p+\frac{3}{4})^{2}+(\frac{k}{2}+\frac{7}{16})}-(\frac{3}{4}+p)$

$=\sqrt{a^{2}+b}-a$,

\noindent
we know that $0\leq\sqrt{1+\frac{n}{2}}-\frac{3}{4}-p<1$ by the above result,
and so $\lfloor\sqrt{1+\frac{n}{2}}-\frac{3}{4}\rfloor=p$.
$\Box$

By this result,
for $n\in (2p^{2}+3p):(2p^{2}+7p+4)$,
Conjecture 1.1 implies that $diam(\mathcal{R}_{n})=n-p$,
and so it means that the diameter of $\mathcal{R}_{n}$ has a liner rule with $n$ in an interval of length $4p+5(=(2p^{2}+7p+4)-(2p^{2}+3p)+1)$ for every $p\geq0$.

Let $p=6$.
Then $2p^{2}+3p+1=91$.
If Conjecture 1.1 is true,
then $diam(\mathcal{R}_{91})=n-p=91-6=85$ by Lemma 4.1.
However from Example 3.7 and Eq. (3.2),
we know that for $n=91$,
$diam(\mathcal{R}_{n})
\geq n-(p+m)+r-2
=91-(6+3)+6-2=86.$
This means that Conjecture 1.1 does not hold for $n=91$.
From Example 3.7,
Eq. (3.2),
and Lemma 4.1,
we can see that if $n-p+r-(2+m)>n-p$,
that is $r>m+2$,
then $diam(\mathcal{R}_{n})>n-p$,
and so Conjecture 1.1 does not hold for $\mathcal{R}_{n}$.
Based on this idea,
we will show that for almost all $n\in\mathbb{N}$,
Conjecture 1.1 does not hold.
First, for a give $p$,
we need to find the ranges of $m$ and $r$ such that $r>m+2$.

\trou \noi {\bf Lemma 4.2.}
\emph{Let $p\geq0$, $m\geq0$ and $n$ be an integer such that $2p^2+3p\leq n\leq2p^2+7p+4$.
Suppose that $\mu$ and $\nu$ are $H$-type strings ending with $p+m$ of length $n$,
and $\alpha$ and $\beta$ be a pair of H-type strings obtained from $\mu$ and $\nu$.}

\noindent
(1)
\emph{If there is a distance-barrier $Bar_{1}$ between $\alpha$ and $\beta$ for $b=m+2+e$ and $e\geq1$,
then}
$\lceil\frac{8p-37-\sqrt{64p^{2}-432p+89}}{10}\rceil\leq m\leq\lfloor\frac{8p-37+\sqrt{64p^{2}-432p+89}}{10}\rfloor$,
\emph{and the maximum value of $e$ is}
\emph{$\lfloor\frac{-6m-11+\sqrt{16m^{2}+(32p+56)m+32p+33}}{6}\rfloor$.}

\noindent
(2)
\emph{If there is a distance-barrier $Bar_{2}$ between $\alpha$ and $\beta$ for $b=m+1+e$ and $e\geq1$,
then}
$\lceil\frac{8p-31-\sqrt{64p^{2}-336p-79}}{10}\rceil\leq m\leq\lfloor\frac{8p-31+\sqrt{64p^{2}-336p-79}}{10}\rfloor$,
\emph{and the maximum value of $e$ is}
\emph{$\lfloor\frac{-6m-9+\sqrt{16m^{2}+(32p+56)m+32p+17}}{6}\rfloor$.}

\trou \noi {\bf Proof.}
First we show that (1) holds.
As there is a
$Bar_{1}=$
$\begin{bmatrix}
(111)^{b-1}1\\
(100)^{b-1}1
\end{bmatrix}$
between $\alpha$ and $\beta$,
the first run of the $H$-type strings $\mu$ (or $\nu$) is a factor of $0^{3b}$ by Lemma 3.8.
Let $\mu$ be such a string with the longest length.
Then we can assume that

$\mu=0^{3b}1^{3b+1}0^{3b+2}\cdots$, and

$\nu=1^{3b}0^{3b+1}1^{3b+2}\cdots$

By the following two facts:

\noindent
$(i)$ the length of a complete $H$-type string ending with $p+m$ is $2(p+m)^{2}+7(p+m)+4$ and the length of
$\mu$ is $n$ such that $2p^2+3p\leq n\leq2p^2+7p+4$;

\noindent
$(ii)$ $\mu$ is the suffix of a complete $H$-type string ending with $p+m$ which is obtained by deleting the factor $10^{2}\cdots 1^{3b-1}$ from this complete $H$-type string,

\noindent
we know that

$(2(p+m)^{2}+7(p+m)+4)-(2p^2+3p)\geq1+2+3+\ldots+(3b-1)$.~~~~~~~~~~~~~~~~~~~~~~~~~~~~$(a)$

Let $b=m+2+e$ and $e\geq1$.
Then

$9e^2+(33+18m)e+(5m^2+19m-8mp-8p+22)\leq 0$~~~~~~~~~~~~~~~~~~~~~~~~~~~~~~~~~~~~~~$(b)$

\noindent
by the above inequality $(a)$.
It is easily seen that the maximum $e$ is the positive root such that the equality holds in inequality $(b)$,
and so $e\leq \frac{-6m-9+\sqrt{16m^{2}+(32p+56)m+32p+17}}{6}$,
as desired.

To get the variation range of $m$,
let $e=1$.
Then by the above inequality $(b)$,
we know that $5m^2+(37- 8p)m+(64-8p)\leq0$,
and so $\frac{8p-37-\sqrt{64p^{2}-432p+89}}{10}\leq m\leq\frac{8p-37+\sqrt{64p^{2}-432p+89}}{10}$,
as desired.

Using the similar method,
we can show that (2) holds.
By Lemma 3.9,
we know that

$(2(p+m)^{2}+7(p+m)+4)-(2p^2+3p)\geq1+2+3+\ldots+3b+(3b+1)+2$.~~~~~~~~~~~~~~~~~~~~~$(c)$

Let $b=m+1+e$,
where $e\geq1$.
Then by inequality $(c)$ we have

$9e^2+(27+18m)e+5m^2+(13-8p)m-8p+16\leq0.$~~~~~~~~~~~~~~~~~~~~~~~~~~~~~~~~~~~~~~$(d)$

From inequality $(d)$,
we know that
$e\leq\frac{-6m-9+\sqrt{16m^{2}+(32p+56)m+32p+17}}{6}$,
as desired.

Let $e=1$.
Then by inequality $(d)$,
we have
$5m^2+(31-8p)m+52-8p\leq0.$
Therefore,
$\frac{8p-31-\sqrt{64p^{2}-336p-79}}{10}\leq m\leq\frac{8p-31+\sqrt{64p^{2}-336p-79}}{10}$,
as desired.
$\Box$

By Eq. (3.2),
we know the value $e$ in Lemma 4.2 is exactly the difference between $d_{\mathcal{R}_{n}}(\alpha,\beta)$ and $n-p$,
that is,
$e=d_{\mathcal{R}_{n}}(\alpha,\beta)-(n-p)=r-(m+2)$,
where $r$ is the number of barriers contained in the distance-barrier between $\alpha$ and $\beta$.

Now we turn to consider Conjecture 1.1.
For convenience,
set

$S_{1}=\{n|n\geq230\}$, and

$S_{2}=\{91,94,95,119:124,131,136,152:170,172,178,181,184,185,189:223\}$.

\trou \noi {\bf Theorem 4.3.}
\emph{Let $n\in S_{1}\cup S_{2}$.
Then $diam(\mathcal{R}_{n})>n-\lfloor\sqrt{1+\frac{n}{2}}-\frac{3}{4}\rfloor$.}

\trou \noi {\bf Proof.}
By Lemma 4.1,
$n-\lfloor\sqrt{1+\frac{n}{2}}-\frac{3}{4}\rfloor=n-p$ if $2p^2+3p\leq n\leq2p^2+7p+4$.
So we need to show that for $n\in S_{1}\cup S_{2}$,
$diam(\mathcal{R}_{n})>n-p$.
The basic idea of the proof is as follows.
First,
for a given $p$ we find $m$ such that there are $H$-type strings $\mu$ and $\nu$ ending with $p+m$ by Lemma 4.2.
Then,
based on $\mu$ and $\nu$,
the $NH$-type strings $\alpha$ and $\beta$ can be constructed by Lemmas 3.8 or 3.9.
Thus,
we know that $diam(\mathcal{R}_{n})\geq d_{\mathcal{R}_{n}}(\alpha,\beta)>n-p$ by Eq. (3.2).

First we show that for {\boldmath$n\in\{91,94,95,131,136,172,178,181,184,185\}$},
the result holds.
For a given $p$,
we need to $H$-type strings $\mu$ and $\nu$ ending with $p+m$ by Lemma 4.2 (2).
By Lemma 3.9,
$NH$-type strings $\alpha$ and $\beta$ can be constructed from $\mu$ and $\nu$,
and for $b=m+2$,
there is a distance-barrier
$\begin{bmatrix}
(111)^{b}1\\
(100)^{b}1
\end{bmatrix}$
between $\alpha$ and $\beta$.
So from $(c)$,
$diam(\mathcal{R}_{n})>n-p$ holds for $n=2(p+m)^{2}+7(p+m)+4-\mathop{\sum}\limits^{3b+1}\limits_{k=1}k-2$.
For $p\leq5$,
we find there is no such $m$.
For $p=6$,
we have $m=1,2$ and $3$.
Therefore $n=94,95$ and $91$.
For $n=94$,
the $H$-type strings (ending with $p+m=6+1=7$) are

$\mu=0^{9}1^{12}0^{13}1^{14}0^{15}1^{16}0^{8}0^{7}$ and

$\nu=1^{9}0^{12}1^{13}0^{14}1^{15}0^{16}1^{8}0^{7}$.

\noindent
By Lemma 3.9,
the $NH$-type strings based on $\mu$ and $\nu$ can be constructed as follows

$\alpha=(100)^{3}10^{11}0^{13}1^{14}0^{15}1^{16}0^{8}0^{7}$ and

$\beta=(111)^{3}10^{11}1^{13}0^{14}1^{15}0^{16}1^{8}0^{7}$.

\noindent
Note that there are $r=4$ barriers contained in the distance-barrier between $\alpha$ and $\beta$ and so $e=r-(m+2)=1$.
So by Eq. (3.2),
$diam(\mathcal{R}_{n})\geq d_{\mathcal{R}_{n}}(\alpha,\beta)=(n-p)+e>n-p=88$.
Similarly,
for $n=95,91$ it can be shown that $e=1$,
and so $diam(\mathcal{R}_{n})\geq d_{\mathcal{R}_{n}}(\alpha,\beta)=(n-p)+e>n-p$.
Note that the case $n=91$ is shown in Example 3.7.
To save space,
for other $p$,
we only give the $m$ and $n$ such that $e\geq1$,
but omit listing the $H$-type strings and the corresponding $NH$--type strings.
For $p=7$,
by Lemma 4.2 (2) we get $m\in0:5$.
If we choose $m=1$ and $2$,
then $n=131$ and $136$,
respectively.
(Note that if $m=0,3,4$ or $5$,
then $n=121,136,131$ or $121$,
respectively.
Obviously,
everyone of those numbers of $n$ is contained in some interval of $S_{1}\cup S_{2}$,
or the same one as the corresponding value of some $m$ that has been list.
The same applies to the other cases.)
For $p=8$,
$m\in0:6$ by Lemma 4.2 (2).
Let $m=1,2,3,4$ and $5$.
Then $n=172,178,181,184$ and $185$.

Next,
we show that the result holds for other integers in $S_{1}\cup S_{2}$.
For a given $p$,
we will first find the $m$ such that $n-p+b-(m+2)>n-p$ by Lemma 4.2 (1),
where $b=m+3$.
Then we know that $n\in(2(p+m)^{2}+7(p+m)+4-\mathop{\sum}\limits^{3b}\limits_{k=1}k-1):
(2(p+m)^{2}+7(p+m)+4-\mathop{\sum}\limits^{3b-1}\limits_{k=1}k)$ by inequality $(a)$,
and so by Lemma 3.8,
the pair of $NH$-type strings can be constructed.
For $p\leq6$,
we find there is not $m$ such that $n-(p+m)+b-2>n-p$.
For $p=7$,
we have $m\in1:3$.
If $m=2$,
then $b=m+3=5$.
For $k\in10:15$ (and so $119\leq n\leq124$) we chose $H$-type strings ending with $p+m=9$

$\mu=0^{k}1^{16}0^{17}1^{18}0^{19}1^{20}0^{10}0^{9}$, and

$\nu=1^{k}0^{16}1^{17}0^{18}1^{19}0^{20}1^{10}0^{9}$.

\noindent
Then by Lemma 3.8,
$NH$-type strings $\alpha_{i}$ and $\beta_{i}$ containing distance-barrier $Bar_{1}$ can be constructed as follows for $k=i+9$,
$i\in1:6$:

$\alpha_{1}=(100)^{4}11^{14}0^{16}1^{18}0^{19}1^{20}0^{10}0^{9}$ and
$\beta_{1}=(111)^{4}10^{14}1^{16}0^{18}1^{19}0^{20}1^{10}0^{9}$,

$\alpha_{2}=(100)^{4}11^{14}0^{17}1^{18}0^{19}1^{20}0^{10}0^{9}$ and
$\beta_{2}=(111)^{4}10^{14}1^{17}0^{18}1^{19}0^{20}1^{10}0^{9}$,

$\alpha_{3}=(100)^{4}11^{15}0^{17}1^{18}0^{19}1^{20}0^{10}0^{9}$ and
$\beta_{3}=(111)^{4}10^{15}1^{17}0^{18}1^{19}0^{20}1^{10}0^{9}$,

$\alpha_{4}=0(100)^{4}11^{15}0^{17}1^{18}0^{19}1^{20}0^{10}0^{9}$ and
$\beta_{4}=1(111)^{4}10^{15}1^{17}0^{18}1^{19}0^{20}1^{10}0^{9}$,

$\alpha_{5}=11000(100)^{3}11^{15}0^{17}1^{18}0^{19}1^{20}0^{10}0^{9}$ and
$\beta_{5}=01111(111)^{3}10^{15}1^{17}0^{18}1^{19}0^{20}1^{10}0^{9}$,

$\alpha_{6}=(100)^{5}1^{16}0^{17}1^{18}0^{19}1^{20}0^{10}0^{9}$ and
$\beta_{6}=(111)^{5}0^{16}1^{17}0^{18}1^{19}0^{20}1^{10}0^{9}$.

\noindent
So we have $b=5$ and so $e=b-(m+2)=1$.
Therefore,
$d_{\mathcal{R}_{n}}(\alpha,\beta)=n-p+e=n-p+1>n-p$ for {\boldmath$n\in119:124$}.

Similarly,
we can consider the cases that $p=8,9$ and $10$.
To save space,
we only list the $H$-type strings which can be obtained by Lemma 4.2(1),
and then the $NH$-type strings which contain distance-barrier $Bar_{1}$ can be constructed by Lemma 3.8.

For $p=8$,
we have $m\in1:3$.
Let $m=1$ and so $b=m+3=4$.
We can chose $H$-type strings ending with $9 (=p+m)$

$\mu=0^{k}1^{13}0^{14}1^{15}0^{16}1^{17}0^{18}1^{19}0^{20}1^{10}0^{9}$ and

$\nu=1^{k}0^{13}1^{14}0^{15}1^{16}0^{17}1^{18}0^{19}1^{20}0^{10}0^{9}$,

\noindent
where $k\in1:12$ (and so $152\leq n\leq163$).
Let $m=3$.
Then $b=m+3=6$ and we can chose $H$-type strings ending with $11(=p+m)$:

$\mu=0^{k}1^{19}0^{20}1^{21}0^{22}1^{23}0^{24}1^{12}0^{11}$ and

$\nu=1^{k}0^{19}1^{20}0^{21}1^{22}0^{23}1^{24}0^{12}0^{11}$,

\noindent
where $k\in12:18$ (and so $164\leq n\leq170$).
So for {\boldmath$n\in152:170$},
we have $diam(\mathcal{R}_{n})\geq n-7>n-8.$

For $p=9$,
we have $m\in0:7$.
Let $m=2$.
Then $b=m+3=5$ and the $H$-type strings can be chosen as ending with $11(=p+m)$:

$\mu_{1}=0^{k_{1}}1^{17}0^{18}1^{19}0^{20}1^{21}0^{22}1^{23}0^{24}1^{12}0^{11}$ and

$\nu_{1}=1^{k_{1}}0^{17}1^{18}0^{19}1^{20}0^{21}1^{22}0^{23}1^{24}0^{12}0^{11}$,

$\mu_{2}=0^{k_{2}}1^{16}0^{17}1^{18}0^{19}1^{20}0^{21}1^{22}0^{23}1^{24}0^{12}0^{11}$ and

$\nu_{2}=1^{k_{2}}0^{16}1^{17}0^{18}1^{19}0^{20}1^{21}0^{22}1^{23}0^{24}1^{12}0^{11}$,

\noindent
where $k_{1}\in2:16$ and $k_{2}\in1:15$ (so {\boldmath$189\leq n\leq203$} and {\boldmath$204\leq n\leq218$}).
Let $m=3$.
Then $b=m+3=6$ and the $H$-type strings can be chosen as ending with $12(=p+m)$:

$\mu=0^{k}1^{19}0^{20}1^{21}0^{22}1^{23}0^{24}1^{25}0^{26}1^{13}0^{12}$ and

$\nu=1^{k}0^{19}1^{20}0^{21}1^{22}0^{23}1^{24}0^{25}1^{26}0^{13}0^{12}$,

\noindent
where $k\in14:18$ (so {\boldmath$219\leq n\leq223$}).

For $p=10$,
we have $0\leq m\leq8$.
Let $m=2$.
Then $b=m+3=5$ and the $H$-type strings can be chosen as ending with $12(=p+m)$:

$\mu_{1}=0^{k_{1}}1^{18}0^{19}1^{20}0^{21}1^{22}0^{23}1^{24}0^{25}1^{26}0^{13}0^{12}$ and

$\nu_{1}=1^{k_{1}}0^{18}1^{19}0^{20}1^{21}0^{22}1^{23}0^{24}1^{25}0^{26}1^{13}0^{12}$,

$\mu_{2}=0^{k_{2}}1^{17}0^{18}1^{19}0^{20}1^{21}0^{22}1^{23}0^{24}1^{25}0^{26}1^{13}0^{12}$ and

$\nu_{2}=1^{k_{2}}0^{17}1^{18}0^{19}1^{20}0^{21}1^{22}0^{23}1^{24}0^{25}1^{26}0^{13}0^{12}$,

$\mu_{3}=0^{k_{3}}1^{16}0^{17}1^{18}0^{19}1^{20}0^{21}1^{22}0^{23}1^{24}0^{25}1^{26}0^{13}0^{12}$ and

$\nu_{3}=1^{k_{3}}0^{16}1^{17}0^{18}1^{19}0^{20}1^{21}0^{22}1^{23}0^{24}1^{25}0^{26}1^{13}0^{12}$,

\noindent
where $k_{1}\in8:18$,
$k_{2}\in1:17$ and
$k_{3}\in1:16$
(so {\boldmath$230\leq n\leq240, 241\leq n\leq256$} and {\boldmath$257\leq n\leq271$}).
Let $m=3$.
Then $b=m+3=6$ and the $H$-type strings can be chosen as ending with $13(=p+m)$:

$\mu=0^{k}1^{19}0^{20}1^{21}0^{22}1^{23}0^{24}1^{25}0^{26}1^{27}0^{28}1^{14}0^{13}$ and

$\nu=1^{k}0^{19}1^{20}0^{21}1^{22}0^{23}1^{24}0^{25}1^{26}0^{27}1^{28}0^{14}0^{13}$,

\noindent
where $k\in9:11$ (so {\boldmath$272\leq n\leq274$}).

For the case $p\geq11$,
we consider the suffix strings $\mu$ and $\nu$ of length $n$ obtained from the complete $H$-type strings ending with $p+2$,
where $2p^{2}+7p+4\geq n \geq 2p^{2}+3p$.
We will show that $NH$-type strings $\alpha$ and $\beta$ can be constructed from $\mu$ and $\nu$ such that there is a distance-barrier
$\begin{bmatrix}
(111)^{5-1}1\\
(100)^{5-1}1
\end{bmatrix}$
between $\alpha$ and $\beta$.
So we have $diam(\mathcal{R}_{n})\geq d_{\mathcal{R}_{n}}(\alpha,\beta)=n-(p+2)+5-2=n-p+1$ by Eq. (3.2).
First,
by the fact that $p\geq11$ and $2(p+2)^2+7(p+2)+4-2(p+1)^2-3(p+1)=8p+21>\mathop{\sum}\limits^{14}\limits_{k=1}k=105$,
we know that if $2p^{2}+7p+4\geq n (\geq 2p^{2}+3p)$,
then the length of the second run of $\mu$ more than $15$.
In other words,
the first two runs of $\mu$ and $\nu$ are $0^{k}1^{3t+s+1}$ and $1^{k}0^{3t+s+1}$,
respectively,
where $k\geq15$ and $s\in0:2$.
By $\mathop{\sum}\limits^{3t+s-1}\limits_{k=1}k<8p+21$,
we know that $t\leq\frac{1-2s+\sqrt{64p+169}}{6}$.
Let $b=5$.
Then $\frac{3b-k-4}{3(t-b)+s+1}\leq\frac{15-1-4}{0+1}=10$.
This means that if we show that the length of the longest run of $\mu$ is at least $3t+s+10$,
then by Lemma 3.8,
we can get a pair of $NH$-type strings $\alpha$ and $\beta$ with the distance-barrier we desire between them.
In fact,
as $\mu$ and $\nu$ end with $p+2$,
the length of the longest run of $\mu$ is $2(p+2+1)=2p+6$.
Thus there is $2p+6-(3t+s+1)$ runs behind the second run (of $\mu$).
Note that

$2p+6-(3t+s+1)$

$\geq 2p+6-s-1-3(\frac{1-2s+\sqrt{64p+169}}{6})$

$=2(p+2)-\frac{1+\sqrt{64p+169}}{2}$.

\noindent
Let
$L(p)=2(p+2)-\frac{1+\sqrt{64p+169}}{2}$.
Then  $L(p)$ is increase with $p$ and $L(11)\geq 11>10$. Thus
we know that for $p\geq11$ and $n\in 2p^{2}+3p:2p^{2}+7p+4$,
the pair of $NH$-type strings $\alpha$ and $\beta$ of length $n$ can be constructed from $H$-type strings $\mu$ and $\nu$ ending with $p+2$,
where there is a distance-barrier
$\begin{bmatrix}
(111)^{5-1}1\\
(100)^{5-1}1
\end{bmatrix}$ between $\alpha$ and $\beta$.
So $diam(\mathcal{R}_{n})\geq d_{\mathcal{R}_{n}}(\alpha,\beta)=n-p+1>n-p$ for {\boldmath$n\geq275$}$=2\times11^{2}+3\times11$.
This completes the proof.
$\Box$

\subsection{Further discussions on $diam(\mathcal{R}_{n})$}

In the above subsection,
we show that Conjecture 1.1 is not true in general meaning by improving the lower bound of the diameter of $\mathcal{R}_{n}$. The tool used is to find  the distance-barriers with the thickness of every barrier being equal to 1.
A more natural idea is to use the distance-barriers with barriers of thickness more than 1 to improve the lower bound of $diam(\mathcal{R}_{n})$ and even to obtain $diam(\mathcal{R}_{n})$.
Some discussions on $diam(\mathcal{R}_{n})$ are given along this line in this subsection.

In the proof of Theorem 4.3,
we show that for $n\geq 275$,
there is always distance-barrier
$\begin{bmatrix}
(111)^{4}1\\
(100)^{4}1
\end{bmatrix}$
between some pairs of $NH$-type strings $\alpha$ and $\beta$ and so we know that $d_{\mathcal{R}_{n}}(\alpha,\beta)=n-p+e$ for $e=1$.
By Lemmas 3.8 and 3.9,
there may exist distance-barriers with more barriers between some pair of vertices $\alpha$ and $\beta$ of $\mathcal{R}_{n}$
and so $d_{\mathcal{R}_{n}}(\alpha,\beta)=n-p+e$ for some $e>1$.

By Lemma 4.2 (1),
for a given $p$ we can find the variation range of $m(\geq0)$ such that there a distance-barrier
$\begin{bmatrix}
(111)^{m+2}1\\
(100)^{m+2}1
\end{bmatrix}$
between the pair of $NH$-type strings $\alpha$ and $\beta$ constructed by Lemma 3.8,
and further for a given $m$,
the maximum $e(\geq1)$ can be obtained such that the distance-barrier
$\begin{bmatrix}
(111)^{m+1+e}1\\
(100)^{m+1+e}1
\end{bmatrix}$ is contained in the pair of $NH$-type strings.
Let $i\in1:e$.
Then the maximum length of the $NH$-type strings $\alpha$ and $\beta$ containing the distance-barrier
$\begin{bmatrix}
(111)^{m+1+i}1\\
(100)^{m+1+i}1
\end{bmatrix}$ is $2(p+m)^{2}+7(p+m)+4-\mathop{\sum}\limits^{3(m+2+i)-1}\limits_{k=1}k$.
Similarly,
by Lemma 4.2 (2),
for a given $p$ we can calculate the $m$ and $e$,
and the length of the $NH$-type strings containing the distance-barrier
$\begin{bmatrix}
(111)^{m+i}1\\
(100)^{m+i}1
\end{bmatrix}$ is $2(p+m)^{2}+7(p+m)+4-\mathop{\sum}\limits^{3(m+1+i)+1}\limits_{k=1}k-2$,
where $i\in1:e$.
Note that $e=d_{\mathcal{R}_{n}}(\alpha,\beta)-(n-p)$.
So the higher this value $e$,
the more lower bound of $diam(\mathcal{R}_{n})$ can be improved.

In the following, first,
for $p=15$ we give an example to improve the lower bound of $diam(\mathcal{R}_{n})$ by finding $e$ with a higher value.
Then by applying the same method, we get the result for $1\leq p\leq 29$.
To save space,
we only list the result.

\trou \noi {\bf Example 4.4.}
\emph{Let $p=15$ and $S=\{518,523,529,532,535,536\}$.
Then $diam(\mathcal{R}_{n})\geq n-13$ for $n\in510:559\setminus S$,
and $diam(\mathcal{R}_{n})\geq n-12$ for $n\in495:509\cup S$.}

\noindent
Note that for $p=15$,
$diam(\mathcal{R}_{n})$ should be $n-15$ by Conjecture 1.1.
The result shown in the above example can be obtained via the following analysis.
By Lemma 4.2 (1),
we calculate $m$, $e$ and the corresponding maximum $n$ as shown in Table 1.
Note that we only consider $2p^{2}+3p(=495)\leq n\leq2p^2+7p+4(=559)$,
and let $n=2p^2+7p+4$ when $n>2p^2+7p+4$ in this process.

\begin{table}[!htbp]
\footnotesize
\caption{$m$ and $e$ when $p=15$ by Lemma 4.2(1)}
\centering
\setlength{\tabcolsep}{1.5mm}{
\begin{tabular}{rcccccccccccccccccc}
\toprule
$m$ &0&1&2&3&4&5&6&7&8&9&10&11&12&13&14&15&16&17\\
\midrule
$e$ &1&2&2&3&3&3&3&3&2&2&2&2&2&1&1&1&1&1\\
\midrule
$n(\leq)$ &523&523&548&502&508&509&505&496&559&559&559&548&523&559&559&559&547&506\\
\bottomrule
\end{tabular}}
\end{table}

\noindent
By Table 1,
we know that $e=3$ for $495\leq n \leq509$ and $e=2$ for $510\leq n\leq559$.
By Lemma 4.2 (2),
we calculate $m$, $e$ and the corresponding $n$ as shown in Table 2.
Note that $495\leq n\leq559$,
when $n>559$ we set $n$ is a non-number.

\begin{table}[!htbp]
\footnotesize
\caption{$m$ and $e$ when $p=15$ by Lemma 4.2(2)}
\centering
\setlength{\tabcolsep}{1.5mm}{
\begin{tabular}{rccccccccccccccccccc}
\toprule
$m$ &0&1&2&3&4&5&6&7&8&9&10&11&12&13&14&15&16&17&18\\
\midrule
$e$ &2&2&3&3&3&3&3&3&3&3&2&2&2&2&2&1&1&1&1\\
\midrule
$n(=)$ &502&535&509&523&532&536&535&529&518&502&-&-&-&541&509&-&-&-&520\\
\bottomrule
\end{tabular}}
\end{table}

The contrast of Tables 1 and 2 shows that we only need to consider $e=3$ for $n>509$ in Table 2.
So we know that $e=3$ for $n=518,523,529,532,535$ and $536$.
This is the end of Example 4.4.

For $p\in 6:29$,
the results are shown in Table 3,
where $d$ denotes the diameter of $\mathcal{R}_{n}$.
Note that for $p\in 6:29$,
$diam(\mathcal{R}_{n})=n-p$ by Conjecture 1.1 if $2p^2+3p\leq n\leq2p^2+7p+4$.
In Table 3,
it is shown that for all $n\in 90:1889$,
the value of $e,e-1$ or $e-2$ are more than $n-p$.

\begin{table}[!htbp]
\footnotesize
\caption{A lower bound of the diameter of $\mathcal{R}_{n}$  when $n\in 90:1889$}
\centering
  \begin{tabular}{|c|c|c|c|}
  \hline
  \diagbox{$p,e$}{$n\in$}{$d\geq$} & $ n-p+e$ & $n-p+e-1$ & $n-p+e-2$\\
  \hline
  6,1 &\{91,93,94\}$(=T_{1})$&90:118$\setminus T_{1}$ & \\
  \hline
  7,1 &119:124$\cup$\{131,136\}$(=T_{2})$&125:151$\setminus T_{2}$ & \\
  \hline
  8,1 &152:170$\cup$\{172,178,181,184,185\}$(=T_{3})$&171:188$\setminus T_{3}$ & \\
  \hline
  9,1 &189:223 &224:229& \\
  \hline
  10,2 &\{230,232,238,239,241\}$(=T_{4})$&230:274$\setminus T_{4}$ & \\
  \hline
  11,2 &275:284$\cup\{230,232,238,239,241\}(=T_{5})$&285:323$\setminus T_{5}$& \\
  \hline
  12,2 &324:352$\cup\{356,361,367,370,373,374\}(=T_{6})$& 353:376$\setminus T_{6}$&\\
  \hline
  13,2 &377:427& 428:433&\\
  \hline
  14,3 &$\{436,446,451\} (=T_{7})$8&434:494 $\setminus T_{7}$&\\
  \hline
  15,3 &495:509$\cup$\{518,523,529,532,535,536\} $(=T_{8})$&510:559$\setminus T_{8}$&\\
  \hline
  16,3 &560:598$\cup\{604,607,617,619,625,626\} (=T_{9})$& 599:628$\setminus T_{9}$&\\
  \hline
  17,3 &629:695m& 696:701&\\
  \hline
  18,4 &$\{703,706,716,718,724,725,727\}(=T_{10})$& 702:778$\setminus T_{10}$&\\
  \hline
  19,4 &779:799$\cup\{806,815,820,826,829,832,833\}(=T_{11})$& 800:859$\setminus T_{11}$&\\
  \hline
  20,4 &860:910$\cup\{916,931,941,945\}(=T_{12})$& 911:944$\setminus T_{12}$&\\
  \hline
  21,5 &$\{946\}(=T_{13})$& 945:1028$\setminus T_{13}$&1029:1033\\
  \hline
  22,5 &$\{1040,1049,1054,1060,1063,1066,1067\}(=T_{14})$& 1034:1126$\setminus T_{14}$& \\
  \hline
  23,5 &1127:1153$\cup\{1157,1171,1174,1192,1193,1195\}(=T_{15})$& 1154:1223$\setminus T_{15}$& \\
  \hline
  24,5 &1224:1286$\cup\{1288,1292,1306,1309,1319,1321\}(=T_{16})$& 1287:1324$\setminus T_{16}$& \\
  \hline
  25,6 &$\{1327,1328,1330\}(=T_{17})$& 1325:1426$\setminus T_{17}$&1427:1429 \\
  \hline
  26,6 &$\{1438,1445,1454,1459,1465,1468,1471,1472\}(=T_{18})$& 1430:1538$\setminus T_{18}$& \\
  \hline
  27,6 &1539:1573$\cup\{1591,1606,1616,1621\}(=T_{19})$& 1574:1651$\setminus T_{19}$& \\
  \hline
  28,6 &1652:1727$\cup\{1732,1744,1751,1760,1765\}(=T_{20})$& 1728:1768$\setminus T_{20}$& \\
  \hline
  29,7 &$\{1771,1774,1777,1778\}(=T_{21})$& 1769:1888$\setminus T_{21}$&1889 \\
  \hline
  \end{tabular}
\end{table}

Let $T=\bigcup\limits^{21}\limits_{i=1}T_{i}$.
Then we have the following result by Table 3.

\trou \noi {\bf Corollary 4.5.}
\emph{Let $1889\geq n\geq90$.
Then}

\begin{center}
$diam(\mathcal{R}_{n})\geq
\begin{cases}
n+1-\lfloor\sqrt{\frac{81}{144}+\frac{5}{18}n}-\frac{13}{12}\rfloor, n\in T;\\
n-\lfloor\sqrt{\frac{81}{144}+\frac{5}{18}n}-\frac{13}{12}\rfloor,n\in 90:1889\setminus T.\\
\end{cases}$
\end{center}

In the following example,
some other types of distance-barriers are constructed to show that the lower bound of the diameter of $\mathcal{R}_{n}$ in Corollary 4.5 (or Table 3) could be further improved for some $n$.

\trou \noi {\bf Example 4.6.}
\emph{Let $n=288$.
Then $diam(\mathcal{R}_{n})\geq n-9$.}

Note that if $n=288$,
then $n=2\times11^{2}+7\times11+4-35$.
By Conjecture 1.1,
$diam(\mathcal{R}_{n})=n-11$.
By Theorem 4.5 (or Table 3),
we know that $diam(\mathcal{R}_{n})\geq n-10$.
Now we show that $diam(\mathcal{R}_{n})\geq n-9$.
Let

$\mu=0^{35}1^{40}0^{41}1^{42}0^{43}1^{44}0^{22}0^{21}$ and

$\nu=1^{35}0^{40}1^{41}0^{42}1^{43}0^{44}1^{22}0^{21}$.

\noindent
Obviously,
$\mu$ and $\nu$ are $H$-type strings of length $n=288$.
Then we can construct the following $NH$-type strings from them:

$\alpha=(11000)^{7}111^{38}0^{41}1^{42}0^{43}1^{44}0^{22}0^{21}$ and

$\beta=(11111)^{7}110^{38}1^{41}0^{42}1^{43}0^{44}1^{22}0^{21}$.

\noindent
It is easily seen that there is a distance-barrier
$\begin{bmatrix}
(11000)^{7}11\\
(11111)^{7}11
\end{bmatrix}$
between $\alpha$ and $\beta$.
The thicknesses of all its barriers are $2$.
By Lemma 3.10,
we have $d_{\mathcal{R}_{n}}(\alpha,\beta)=d_{\mathcal{R}_{n}}(\mu,\nu)+2(8-2)=H(\mu\nu)+12=n-21+12=n-9$.
This is the end of Example 4.6.

\trou \noi {\bf Example 4.7.}
\emph{Let $n=520$.
Then $diam(\mathcal{R}_{n})\geq n-12$.}

By Conjecture 1.1,
we have $diam(\mathcal{R}_{n})=d_{\mathcal{R}_{n}}(\mu,\nu)=n-15$ since $n=520=2\times15^{2}+7\times15+4-39$.
By Corollary 4.5 (or Table 3),
we know that $diam(\mathcal{R}_{n})\geq n-13$.
Let

$\mu=0^{38}1^{39}0^{40}1^{41}0^{42}1^{43}0^{44}1^{45}0^{46}1^{47}0^{48}1^{24}0^{23}$ and

$\nu=1^{38}0^{39}1^{40}0^{41}1^{42}0^{43}1^{44}0^{45}1^{46}0^{47}1^{48}0^{24}0^{23}$.

\noindent
Obviously,
$\mu$ and $\nu$ are $H$-type strings of length $n=520$.
Then we can construct the following $NH$-type strings from them:

$\alpha=100(11000)^{7}1^{39}0^{40}1^{41}0^{42}1^{43}0^{44}1^{45}0^{46}1^{47}0^{48}1^{24}0^{23}$ and

$\beta=111(11111)^{7}0^{39}1^{40}0^{41}1^{42}0^{43}1^{44}0^{45}1^{46}0^{47}1^{48}0^{24}0^{23}$.

\noindent
It is easily seen that there is a distance-barrier
$\begin{bmatrix}
100(11000)^{6}11\\
111(11111)^{6}11
\end{bmatrix}$
between $\alpha$ and $\beta$.
The thicknesses of all its barriers are $2$ except the first one is 1.
By Lemma 3.10,
we have $d_{\mathcal{R}_{n}}(\alpha,\beta)
=d_{\mathcal{R}_{n}}(\mu,\nu)+(2(7-2)+1)
=H(\mu,\nu)+12
=n-23+11
=n-12.$

Those examples show that we need to consider all the possible distance-barriers between vertices of $\mathcal{R}_{n}$ to find the diameter of $\mathcal{R}_{n}$.
In fact,
we can find a lower bound of $diam(\mathcal{R}_{n})$ for a given distance-barrier using the method as given in Lemma 4.2.
However,
there are many types of distance-barriers which are effective to the (maximum) distance of $\mathcal{R}_{n}$,
as shown in Examples 4.6 and 4.7.
We believe it is not an easy task to find such distance-barriers.

Being inclined to the view that there is no unified formula of $diam(\mathcal{R}_{n})$, we end this section with the following open problem.

\trou \noi {\bf Question 4.8.} \emph{What is the diameter of $\mathcal{R}_{n}$?}

\end{document}